\documentclass[12pt,reqno]{amsart}
\usepackage{amsmath,amssymb,amscd,latexsym,amsthm,mathrsfs,comment}
\usepackage[unicode]{hyperref}
\usepackage{hypbmsec}
\ifx\pdfoutput\undefined\else
\hypersetup{
colorlinks=true,
linkcolor=blue,
citecolor=blue,
urlcolor=blue,
filecolor=blue,
bookmarksnumbered=true,
pdfstartview=FitH,
pdfhighlight=/N
}
\fi
\newtheorem{theorem}{Theorem}
\newtheorem{lemma}{Lemma}

\theoremstyle{remark}
\newtheorem{remark}{Remark}

\renewcommand\pmod[1]{\;(\operatorname{mod}#1)}
\newcommand\ord{\operatorname{ord}}

\begin{document}

\hypersetup{pdfauthor={Wadim Zudilin},%
pdftitle={Ramanujan-type supercongruences}}

\title{Ramanujan-type supercongruences}

\author{Wadim Zudilin}
\address{Max-Planck-Institut f\"ur Mathematik,
Vivatsgasse~7, D-53111 Bonn, GERMANY\newline\hbox to\parindent{\hss}%
\textit{WWW-address}: \href{http://wain.mi.ras.ru/}{\tt http://wain.mi.ras.ru/}}

\date{October 14, 2008}

\thanks{The work was supported by a fellowship
of the Max Planck Institute for Mathematics (Bonn)
and by project no.~R-146-000-103-112 of the
National University of Singapore Academic Research Fund.}

\begin{abstract}
We present several supercongruences that may be viewed as $p$-adic analogues
of Ramanujan-type series for $1/\pi$ and $1/\pi^2$, and prove three of
these examples.
\end{abstract}

\maketitle

\section{Introduction}
\label{s1}

Throughout the paper, the letter~$p$ is used for primes.
It was observed by L.~Van Hamme \cite{VH1} that several Ramanujan's and Ramanujan-like
formulas~\cite{Ra} for $1/\pi$,
\begin{align}
\sum_{n=0}^\infty\frac{(\frac12)_n^3}{n!^3}(4n+1)(-1)^n
&=\frac2\pi,
\label{A01a}
\displaybreak[2]\\
\sum_{n=0}^\infty\frac{(\frac12)_n^3}{n!^3}(6n+1)\frac1{4^n}
&=\frac4\pi,
\nonumber
\displaybreak[2]\\
\sum_{n=0}^\infty\frac{(\frac12)_n^3}{n!^3}(6n+1)\frac{(-1)^n}{8^n}
&=\frac{2\sqrt2}\pi,
\nonumber
\displaybreak[2]\\
\sum_{n=0}^\infty\frac{(\frac12)_n^3}{n!^3}(42n+5)\frac1{64^n}
&=\frac{16}\pi,
\nonumber
\end{align}
admit very nice conjectural $p$-analogues:
\begin{align}
\sum_{n=0}^{p-1}\frac{(\frac12)_n^3}{n!^3}(4n+1)(-1)^n
&\overset?\equiv\biggl(\frac{-1}p\biggr)p\pmod{p^3} \quad\text{for $p>2$},
\label{A01b}
\displaybreak[2]\\
\sum_{n=0}^{p-1}\frac{(\frac12)_n^3}{n!^3}(6n+1)\frac1{4^n}
&\overset?\equiv\biggl(\frac{-1}p\biggr)p\pmod{p^3} \quad\text{for $p>2$},
\label{A02b}
\displaybreak[2]\\
\sum_{n=0}^{p-1}\frac{(\frac12)_n^3}{n!^3}(6n+1)\frac{(-1)^n}{8^n}
&\overset?\equiv\biggl(\frac{-2}p\biggr)p\pmod{p^3} \quad\text{for $p>2$},
\label{A03b}
\displaybreak[2]\\
\sum_{n=0}^{p-1}\frac{(\frac12)_n^3}{n!^3}(42n+5)\frac1{64^n}
&\overset?\equiv5\biggl(\frac{-1}p\biggr)p\pmod{p^3} \quad\text{for $p>2$},
\label{A04b}
\end{align}
where $\bigl(\frac{\cdot}p\bigr)$ and $(a)_n=\Gamma(a+n)/\Gamma(n)$
denote the Legendre (quadratic residue) symbol and the Pochhammer symbol, respectively.
Note that since the $p$-adic order of $(\frac12)_n/n!$
is~1 for $n=\frac{p+1}2,\dots,\allowbreak p-\nobreak1$,
one may compute the sums in~\eqref{A01b}--\eqref{A04b} for $n$ up to~$\frac{p-1}2$
(and this is how the four supercongruences%
\footnote{A congruence modulo a prime~$p$ is said to be a supercongruence
if it happens to hold modulo some higher power of~$p$.}
\eqref{A01b}--\eqref{A04b}
were conjectured in~\cite{VH1}); for example,
an equivalent form of~\eqref{A01b} is~\cite[(B.2)]{VH1}
\begin{equation}
\sum_{n=0}^{(p-1)/2}\frac{(\frac12)_n^3}{n!^3}(4n+1)(-1)^n
\overset?\equiv(-1)^{(p-1)/2}p\pmod{p^3} \quad\text{for $p>2$}.
\label{A01}
\end{equation}
Van Hamme himself was able to prove \eqref{A01b} modulo~$p$ (not~$p^3$);
recently E.~Mortenson~\cite{Mo} gave a proof for~\eqref{A01b}.
As far as we know, no progress has been done towards proving the other
congruences \eqref{A02b}, \eqref{A03b}, and \eqref{A04b}.

\begin{remark}
\label{r1}
In fact, one can state all these $p$-supercongruences as infinite sums
by replacing $n!$ in the denominators of the summands by $(-1)^n\Gamma_p(n+1)$,
where $\Gamma_p$ is Morita's $p$-adic gamma function~\cite{Ma}.
For instance, the congruence~\eqref{A01b} takes the form
\begin{equation}
\sum_{n=0}^\infty\frac{(\frac12)_n^3}{\Gamma_p(n+1)^3}(4n+1)
\equiv(-1)^{(p-1)/2}p\pmod{p^3} \quad\text{for $p>2$}.
\label{A01g}
\end{equation}
\end{remark}

It seems that \emph{all} known Ramanujan-type formulas for $1/\pi$
admit similar $p$-anal\-ogues. One example, which is not
given in the list~\cite{VH1}, is an analogue of the formula
\begin{equation}
\sum_{n=0}^\infty\frac{(\frac12)_n(\frac14)_n(\frac34)_n}{n!^3}(20n+3)\frac{(-1)^n}{2^{2n}}
=\frac8\pi;
\label{x01a}
\end{equation}
it reads as
\begin{equation}
\sum_{n=0}^{p-1}\frac{(\frac12)_n(\frac14)_n(\frac34)_n}{n!^3}(20n+3)\frac{(-1)^n}{2^{2n}}
\overset?\equiv3\biggl(\frac{-1}p\biggr)p\pmod{p^3} \quad\text{for $p>2$}.
\label{x01b}
\end{equation}
Even more, some recently discovered, by J.~Guillera, formulas for $1/\pi^2$
(see \cite{G1}--\cite{G3} and~\cite{Zu}), like
\begin{equation}
\sum_{n=0}^\infty\frac{(\frac12)_n^3(\frac14)_n(\frac34)_n}{n!^5}(120n^2+34n+3)\frac1{2^{4n}}
=\frac{32}{\pi^2},
\label{B03a}
\end{equation}
possess such analogues as well:
\begin{equation}
\sum_{n=0}^{p-1}\frac{(\frac12)_n^3(\frac14)_n(\frac34)_n}{n!^5}(120n^2+34n+3)\frac1{2^{4n}}
\overset?\equiv 3p^2\pmod{p^5} \quad\text{for $p>2$}.
\label{B03b}
\end{equation}
In Section~\ref{s3} we present a (far from exhaustive) list
of many other supercongruences which are analogues of the examples
indicated in our recent mini-survey~\cite{Zu}; there one can find
the analogue~\eqref{B08b} of the formula for~$1/\pi^3$ experimentally
discovered by B.~Gourevich (see~\cite{Zu}). Unfortunately, all these
congruences remain conjectures, except the already mentioned example
\eqref{A01b} and, as we show below, the congruences \eqref{x01b} and~\eqref{B03b}.
We do not dispose of a suitable general method to prove our
experimentally discovered supercongruences. But at least in the three cases
we can use the WZ-method \cite{EZ}, \cite{G1}, \cite{G3} which results
in very simple proofs given in Section~\ref{s2} of the following results.

\begin{theorem}
\label{t1}
The congruence in~\eqref{A01} is valid.
\end{theorem}

Theorem~\ref{t1} is the main result in~\cite{Mo}, but our proof
is essentially different from the one in~\cite{Mo} and shorter.

\begin{theorem}
\label{t2}
The following congruence is true:
\begin{equation}
\sum_{n=0}^{p-1}\frac{(\frac12)_n^3(\frac12)_{2n}}{n!^5}(120n^2+34n+3)\frac1{2^{6n}}
\equiv3p^2\pmod{p^5} \quad\text{for $p>2$};
\label{B03}
\end{equation}
it is equivalent to~\eqref{B03b}.
\end{theorem}

\begin{theorem}
\label{t3}
We have
\begin{equation}
\sum_{n=0}^{p-1}\frac{(\frac12)_n(\frac12)_{2n}}{n!^3}(20n+3)\frac1{2^{4n}}
\equiv3(-1)^{(p-1)/2}p\pmod{p^3} \quad\text{for $p>2$},
\label{x01}
\end{equation}
which is equivalent to~\eqref{x01b}.
\end{theorem}

\section{Proofs by the WZ-machinery}
\label{s2}

The original argument of D.~Zeilberger and S.~Ekhad~\cite{EZ},
applied there to proving the identity~\eqref{A01a} (which, in fact,
starts its history from Bauer's work~\cite{Bau} in~1859, much earlier than
Ramunujan's birth), was later developed in Guillera's works~\cite{G1},~\cite{G3}.
It is based on introducing suitable WZ-pairs and creative telescoping.
This method has rather strong limitations (and it always requires a
preliminary human guess) but it is the only one known so far for proving
the hypergeometric formulas for $1/\pi^2$, like~\eqref{B03a}.
In our proof of Theorems~\ref{t1}, \ref{t2}, and~\ref{t3} we do not need
an intelligent guess, since we can simply borrow
the Zeilberger--Ekhad--Guillera WZ-pairs for~\eqref{A01a}, \eqref{B03a},
and~\eqref{x01a}, respectively.

\begin{proof}[Proof of Theorem~{\rm\ref{t1}}]
Introduce the rational functions
\begin{gather*}
F(n,k)
=(-1)^{n+k}(4n+1)\frac{(\frac12)_n^2(\frac12)_{n+k}}{(1)_n^2(1)_{n-k}(\frac12)_k^2},
\\
G(n,k)
=-\frac{(2n-1)^2}{2(n-k)(4n-3)}F(n-1,k)
=(-1)^{n+k}\cdot2\cdot\frac{(\frac12)_n^2(\frac12)_{n+k-1}}{(1)_{n-1}^2(1)_{n-k}(\frac12)_k^2}
\end{gather*}
in two parameters $n$ and~$k$.
They form a WZ-pair, namely, they satisfy
\begin{equation}
F(n,k-1)-F(n,k)=G(n+1,k)-G(n,k).
\label{00}
\end{equation}
Indeed,
\begin{gather*}
\frac{F(n,k-1)}{F(n,k)}=-\frac{(k-\frac12)^2}{(n+k-\frac12)(n-k+1)},
\quad
\frac{F(n,k)}{G(n,k)}=\frac{(4n+1)(n+k-\frac12)}{2n^2},
\\
\frac{G(n,k+1)}{G(n,k)}=-\frac{(n+\frac12)^2(n+k-\frac12)}{n^2(n-k+1)},
\end{gather*}
and it is routine to verify the identity
$$
-\frac{(k-\frac12)^2(4n+1)}{2(n-k+1)n^2}
-\frac{(4n+1)(n+k-\frac12)}{2n^2}
=-\frac{(n+\frac12)^2(n+k-\frac12)}{n^2(n-k+1)}
-1
$$
which is the result of division of both sides of~\eqref{00} by $G(n,k)$.

Summing \eqref{00} over $n=0,1,\dots,\frac{p-1}2$, we obtain
\begin{equation}
\sum_{n=0}^{(p-1)/2}F(n,k-1)
-\sum_{n=0}^{(p-1)/2}F(n,k)
=G(\tfrac{p+1}2,k)-G(0,k)
=G(\tfrac{p+1}2,k).
\label{01}
\end{equation}
Furthermore, for $k=1,2,\dots,\frac{p-1}2$ we have
\begin{align*}
G(\tfrac{p+1}2,k)
&=(-1)^{(p+1)/2+k}\cdot2\cdot\biggl(\frac{(\frac12)_{(p-1)/2}}{(1)_{(p-1)/2}}\biggr)^2
\frac{(\frac{p-1}2+\frac12)^2(\frac12)_{(p+1)/2+k-1}}{(1)_{(p+1)/2-k}(\frac12)_k^2}
\\
&=(-1)^{(p+1)/2+k}\cdot2^{-p}{\binom{p-1}{\frac{p-1}2}}^2p^2
\cdot\frac{(\frac12)_{(p+1)/2+k-1}}{(1)_{(p+1)/2-k}(\frac12)_k^2}
\equiv 0\pmod{p^3},
\end{align*}
since $(\frac12)_{(p+1)/2+k-1}$ is divisible by~$(\frac12)_{(p+1)/2}$, hence by~$p$,
while the denominator is coprime to~$p$. Comparing this result with~\eqref{01} we see that
\begin{equation}
\sum_{n=0}^{(p-1)/2}F(n,0)
\equiv\sum_{n=0}^{(p-1)/2}F(n,1)
\equiv\sum_{n=0}^{(p-1)/2}F(n,2)
\equiv\dots
\equiv\sum_{n=0}^{(p-1)/2}F(n,\tfrac{p-1}2)\pmod{p^3}.
\label{02}
\end{equation}
On the other hand,
\begin{align}
\sum_{n=0}^{(p-1)/2}F(n,\tfrac{p-1}2)
&=F(\tfrac{p-1}2,\tfrac{p-1}2)
=(4\cdot\tfrac{p-1}2+1)\frac{(\frac12)_{p-1}}{(1)_{(p-1)/2}^2}
\nonumber\\
&=2\frac{(\frac12)_p}{(1)_{(p-1)/2}^2}
=p\cdot 2^{-2(p-1)}\binom{2p-1}{p-1}\binom{p-1}{\frac{p-1}2}.
\label{03}
\end{align}
It remains to use the well-known congruences
\begin{equation*}
\binom{2p-1}{p-1}
\equiv1\pmod{p^3}
\end{equation*}
due to J.~Wolstenholme (see \cite{Wo}) and
\begin{equation}
\binom{p-1}{\frac{p-1}2}
\equiv(-1)^{(p-1)/2}2^{2(p-1)}\pmod{p^3}
\label{03b}
\end{equation}
due to F.~Morley~\cite{My},
although they are true modulo~$p^2$ only if $p=3$,
to conclude that the expression in~\eqref{03} is congruent to
$$
\sum_{n=0}^{(p-1)/2}F(n,\tfrac{p-1}2)
\equiv(-1)^{(p-1)/2}p\pmod{p^4}
$$
(for our purposes we need the latter modulo~$p^3$),
and the desired result follows from~\eqref{02}.
\end{proof}

\begin{proof}[Proof of Theorem~{\rm\ref{t2}}]
The validity of~\eqref{B03} for $p=3$ is checked by hand; therefore,
we may assume that $p\ge5$.

The rational functions
\begin{gather*}
F(n,k)
=(120n^2-84nk+34n-10k+3)\frac{(\frac12)_n^3(\frac12)_{2n+k}}{2^{6n}(1)_n^3(1)_{n-k}^2(\frac12)_k^3},
\\
G(n,k)
=\frac{(\frac12)_n^3(\frac12)_{2n+k-1}}{2^{6n-8}(1)_{n-1}^3(1)_{n-k}^2(\frac12)_k^3}
\end{gather*}
satisfy~\eqref{00}, since after division of both sides by $G(n,k)$ one
only needs to verify that
\begin{align*}
&
\frac{(120n^2-84nk+118n-10k+13)(k-\frac12)^3}{256n^3(n-k+1)^2}
\\ &\qquad
-\frac{(120n^2-84nk+34n-10k+3)(2n+k-\frac12)}{256n^3}
\\ &\quad
=\frac{(2n+k-\frac12)(2n+k+\frac12)(n+\frac12)^3}{64n^3(n-k+1)^2}-1.
\end{align*}
Summing \eqref{00}, this time over $n=0,1,\dots,p-1$, we obtain
$$
\sum_{n=0}^{p-1}F(n,k-1)-\sum_{n=0}^{p-1}F(n,k)
=G(p,k)-G(0,k)=G(p,k).
$$
Using
$$
G(p,k)
=\frac{(\frac12)_p^3(\frac12)_{2p+k-1}}{2^{6n-8}(1)_{p-1}^3(1)_{p-k}^2(\frac12)_k^3}
\equiv0\pmod{p^5}
$$
for $k=1,2,\dots,\frac{p-1}2$, we deduce as above that
\begin{equation*}
\sum_{n=0}^{p-1}F(n,0)
\equiv\sum_{n=0}^{p-1}F(n,\tfrac{p-1}2)\pmod{p^5}.
\end{equation*}
In addition, for $n=\frac{p+1}2,\frac{p+1}2+1,\dots,p-1$ we have
\begin{align*}
\ord_p\frac{(\frac12)_n^3(\frac12)_{2n+(p-1)/2}}{2^{6n}(1)_n^3(1)_{n-(p-1)/2}^2(\frac12)_{(p-1)/2}^3}
&=\ord_p(\tfrac12)_n^3(\tfrac12)_{2n+(p-1)/2}
\\
&\ge\ord_p(\tfrac12)_{(p+1)/2}^3(\tfrac12)_{p+1+(p-1)/2}
=5;
\end{align*}
thus
$$
F(n,\tfrac{p-1}2)\equiv0\pmod{p^5}
$$
for these values of $n$, and we obtain
\begin{equation*}
\sum_{n=0}^{p-1}F(n,\tfrac{p-1}2)
=\sum_{n=(p-1)/2}^{p-1}F(n,\tfrac{p-1}2)
\equiv F(\tfrac{p-1}2,\tfrac{p-1}2)\pmod{p^5}.
\end{equation*}
Finally, note that
$$
F(\tfrac{p-1}2,\tfrac{p-1}2)
=3p(3p-2)\frac{(\frac12)_{3(p-1)/2}}{2^{3(p-1)}(1)_{(p-1)/2}^3}
=3p\cdot 2^{-6(p-1)}\frac{(3p-2)!}{(\frac{3(p-1)}2)!\,(\frac{p-1}2)!^3}
$$
and
\begin{equation}
2^{-6(p-1)}\frac{(3p-2)!}{(\frac{3(p-1)}2)!\,(\frac{p-1}2)!^3}
\equiv p\pmod{p^4},
\label{06}
\end{equation}
the congruence following from the lemma below,
to arrive at the desired claim~\eqref{B03}.
\end{proof}

\begin{lemma}
For a prime $p\ge5$ take $n=(p-1)/2$. Then
\begin{equation*}
\frac{(6n+1)!}{(3n)!\,n!^3}
\equiv p\cdot2^{12n}\pmod{p^4}.
\end{equation*}
\end{lemma}

\begin{proof}
For $p>3$, we have
\begin{align*}
\frac{(6n+1)!}{(3n)!}
&=\frac{(6n+2)!}{2\cdot(3n+1)!}
=(2n+1)\prod_{j=1}^n(2p-j)(2p+j)(3p-j)
\\
&=pn!^3\prod_{j=1}^n\biggl(1-\frac{2p}j\biggr)\biggl(1+\frac{2p}j\biggr)\biggl(1-\frac{3p}j\biggr)
\\
&=pn!^3\Bigl(1-3pH_n+\frac92p^2H_n^2-\frac{17}2p^2H_n'+O(p^3)\Bigr),
\end{align*}
where
$$
H_n=\sum_{l=1}^n\frac1l \quad\text{and}\quad
H_n'=\sum_{l=1}^n\frac1{l^2}\equiv0\pmod p.
$$
Thus,
\begin{equation}
\frac{(6n+1)!}{(3n)!\,n!^3}
\equiv p\Bigl(1-3pH_n+\frac92p^2H_n^2\Bigr)\pmod{p^4}.
\label{06b}
\end{equation}
It remains to use a result of E.~Lehmer~\cite[Eq.~(45)]{Le},
\begin{equation*}
H_n\equiv -2Q_p(2)+pQ_p(2)^2\pmod{p^2},
\qquad\text{where}\quad
Q_p(2)=\frac{2^{p-1}-1}p,
\end{equation*}
which may be written in the form
$$
Q_p(2)\equiv-\frac12H_n+\frac18pH_n^2\pmod{p^2},
$$
equivalently,
\begin{equation*}
2^{p-1}\equiv1-\frac12pH_n+\frac18p^2H_n^2\pmod{p^3}.
\end{equation*}
Raising both sides of the congruence to the sixth power we obtain
$$
2^{12n}=2^{6(p-1)}
\equiv1-3pH_n+\frac92p^2H_n^2\pmod{p^3}.
$$
The comparison of this with \eqref{06b} leads to the desired claim.
\end{proof}

\begin{proof}[Proof of Theorem~{\rm\ref{t3}}]
This time the choice of a WZ-pair satysfying~\eqref{00} is as follows:
\begin{gather*}
F(n,k)
=(20n-2k+3)\frac{(-1)^{n+k}(\frac12)_n(\frac12)_{2n+k}}{2^{4n}(1)_n^2(1)_{n-k}(\frac12)_k^2},
\\
G(n,k)
=\frac{(-1)^{n+k}(\frac12)_n(\frac12)_{2n+k-1}}{2^{4n-6}(1)_{n-1}^2(1)_{n-k}(\frac12)_k^2}.
\end{gather*}
Since $G(p,k)\equiv0\pmod{p^3}$ for $k=1,2,\dots,\frac{p-1}2$, by applying~\eqref{00}
we conclude that
$$
\sum_{n=0}^{p-1}F(n,0)
\equiv\sum_{n=0}^{p-1}F(n,\tfrac{p-1}2)\pmod{p^3}.
$$
Furthermore, $F(n,\frac{p-1}2)\equiv0\pmod{p^3}$ for $n=\frac{p+1}2,\dots,p-1$,
hence the latter sum reduces to the single term $F(\frac{p-1}2,\frac{p-1}2)\pmod{p^3}$.
Finally,
\begin{align*}
F(\tfrac{p-1}2,\tfrac{p-1}2)
&=3(3p-2)\frac{(\frac12)_{3(p-1)/2}}{2^{2(p-1)}(1)_{(p-1)/2}(\frac12)_{(p-1)/2}}
\\
&=3\cdot2^{-6(p-1)}\frac{(3p-2)!}{(\frac{3(p-1)}2)!\,(\frac{p-1}2)!^3}
\cdot2^{2(p-1)}{\binom{p-1}{\frac{p-1}2}}^{-1}
\end{align*}
and application of~\eqref{03b} and~\eqref{06} finishes the proof of~\eqref{x01}.
\end{proof}

\section{Experimental database and conclusion}
\label{s3}

To convince the reader that there is a general pattern for
$p$-supercongruences analogous to Ramanujan-type series for $1/\pi$
and their generalizations, we provide a list of several
examples (the corresponding original formulas, some of them
conjectural as well, are all represented in~\cite{Zu}):
\begin{align}
\sum_{n=0}^{p-1}\frac{(\frac12)_n(\frac16)_n(\frac56)_n}{n!^3}(5418n+263)\frac{(-1)^n}{80^{3n}}
&\overset?\equiv263\biggl(\frac{-15}p\biggr)p\pmod{p^3} \quad\text{for $p>5$},
\label{A05b}
\displaybreak[2]\\
\sum_{n=0}^{p-1}\frac{(\frac12)_n(\frac14)_n(\frac34)_n}{n!^3}(21460n+1123)\frac{(-1)^n}{882^{2n}}
&\overset?\equiv1123\biggl(\frac{-1}p\biggr)p\pmod{p^3} \quad\text{for $p>7$},
\label{A06b}
\displaybreak[2]\\
\sum_{n=0}^{p-1}\frac{(\frac12)_n(\frac14)_n(\frac34)_n}{n!^3}(26390n+1103)\frac1{99^{4n}}
&\overset?\equiv1103\biggl(\frac{-2}p\biggr)p\pmod{p^3} \quad\text{for $p>11$},
\label{A07b}
\displaybreak[2]\\
\sum_{n=0}^{p-1}\frac{(\frac12)_n(\frac13)_n(\frac23)_n}{n!^3}(14151n+827)\frac{(-1)^n}{500^{2n}}
&\overset?\equiv827\biggl(\frac{-3}p\biggr)p\pmod{p^3} \quad\text{for $p>5$},
\label{A08b}
\displaybreak[2]\\
\sum_{n=0}^{p-1}\frac{(\frac12)_n^5}{n!^5}(20n^2+8n+1)\frac{(-1)^n}{2^{2n}}
&\overset?\equiv p^2\pmod{p^5} \quad\text{for $p>3$},
\label{B01b}
\displaybreak[2]\\
\sum_{n=0}^{p-1}\frac{(\frac12)_n^5}{n!^5}(820n^2+180n+13)\frac{(-1)^n}{2^{10n}}
&\overset?\equiv13p^2\pmod{p^5} \quad\text{for $p>3$},
\label{B02b}
\end{align}
\\[-10mm]
\begin{align}
\sum_{n=0}^{p-1}\frac{(\frac12)_n(\frac14)_n(\frac34)_n(\frac16)_n(\frac56)_n}{n!^5}
(1640n^2+278n+15)\frac{(-1)^n}{2^{10n}}
&\overset?\equiv15\biggl(\frac3p\biggr)p^2\pmod{p^5}
\label{B04b}
\\[-6mm] \nonumber &
\quad\text{for $p>3$},
\displaybreak[2]\\
\sum_{n=0}^{p-1}\frac{(\frac12)_n(\frac14)_n(\frac34)_n(\frac13)_n(\frac23)_n}{n!^5}
(252n^2+63n+5)\frac{(-1)^n}{48^n}
&\overset?\equiv 5p^2\pmod{p^5}
\label{B05b}
\\[-2mm] \nonumber &
\quad\text{for $p>3$},
\displaybreak[2]\\
\sum_{n=0}^{p-1}\frac{(\frac12)_n(\frac13)_n(\frac23)_n(\frac16)_n(\frac56)_n}{n!^5}
(5418n^2+693n+29)\frac{(-1)^n}{80^{3n}}
&\overset?\equiv29\biggl(\frac5p\biggr)p^2\pmod{p^5}
\label{B06b}
\\[-6mm] \nonumber &
\quad\text{for $p>5$},
\displaybreak[2]\\
\sum_{n=0}^{p-1}\frac{(\frac12)_n(\frac18)_n(\frac38)_n(\frac58)_n(\frac78)_n}{n!^5}
(1920n^2+304n+15)\frac1{7^{4n}}
&\overset?\equiv15\biggl(\frac7p\biggr)p^2\pmod{p^5}
\label{B07b}
\\[-6mm] \nonumber &
\quad\text{for $p>7$},
\displaybreak[2]\\
\sum_{n=0}^{p-1}\frac{(\frac12)_n^7}{n!^7}
(168n^3+76n^2+14n+1)\frac1{2^{6n}}
&\overset?\equiv\biggl(\frac{-1}p\biggr)p^3\pmod{p^7}
\label{B08b}
\\[-2mm] \nonumber &
\quad\text{for $p>2$},
\displaybreak[2]\\
\sum_{n=0}^{p-1}\sum_{k=0}^n{\binom nk}^4\cdot(4n+1)\frac1{6^{2n}}
&\overset?\equiv\biggl(\frac{-15}p\biggr)p\pmod{p^2}
\label{C01b}
\\[-2mm] \nonumber &
\quad\text{for $p>3$},
\displaybreak[2]\\
\sum_{n=0}^{p-1}\sum_{k=0}^{[n/3]}(-1)^{n-k}3^{n-3k}
\frac{(3k)!}{k!^3}\binom n{3k}\binom{n+k}k
\cdot(4n+1)\frac1{81^n}
&\overset?\equiv\biggl(\frac{-3}p\biggr)p\pmod{p^3}
\label{C02b}
\\[-2mm] \nonumber &
\quad\text{for $p>3$},
\displaybreak[2]\\
\sum_{n=0}^{p-1}\sum_{k=0}^n\binom nk^2\binom{2k}k\binom{2n-2k}{n-k}
\cdot(5n+1)\frac1{64^n}
&\overset?\equiv\biggl(\frac{-3}p\biggr)p\pmod{p^3}
\label{C03b}
\\[-2mm] \nonumber &
\quad\text{for $p>2$},
\displaybreak[2]\\
\sum_{n=0}^{p-1}\binom{2n}n^2\sum_{k=0}^n\binom{2k}k^2\binom{2n-2k}{n-k}^2
\cdot(36n^2+12n+1)\frac1{2^{10n}}
&\overset?\equiv p^2\pmod{p^3}
\label{C04b}
\\[-2mm] \nonumber &
\quad\text{for $p>2$}.
\end{align}

\begin{remark}
\label{r2}
Replacing the finite sums on~$n$ in~\eqref{C01b}--\eqref{C04b}
by infinite sums following the `recipe' in Remark~\ref{r1}
(that is, replacing some factorials by the corresponding
values of the $p$-adic gamma function as in~\eqref{A01g})
looks like a tricky problem. It seems that this is the place
where we lose some `natural' powers of~$p$ in~\eqref{C01b} and~\eqref{C04b}.
\end{remark}

We find it quite unfortunate that, in the other cases when we have
proofs of the original Ramanujan-type series for $1/\pi$ and $1/\pi^2$
by creative telescoping, we cannot recover our argument in Section~\ref{s2}.
These are the supercongruences \eqref{A02b}--\eqref{A04b}, \eqref{B01b}
and~\eqref{B02b}; the choices of the first member in the corresponding
WZ-pair (again borrowed from~\cite{G1},~\cite{G3}) are
\begin{align*}
F(n,k)
&=(6n-2k+1)
\frac{(\frac12)_n(\frac12)_{n+k}(\frac12)_{n-k}}
{2^{2n}(1)_n^2(1)_{n-k}(\frac12)_k},
\displaybreak[2]\\
F(n,k)
&=(6n-2k+1)
\frac{(-1)^{n+k}(\frac12)_{n+k}(\frac12)_{n-k}^2}
{2^{3n-k}(1)_n^2(1)_{n-k}},
\displaybreak[2]\\
F(n,k)
&=(84n^2-56nk+4k^2+52n-12k+5)
\frac{(-1)^k(\frac12)_n(\frac12)_{n+k}(\frac12)_{n-k}^2}
{2^{4n}(1)_n^2(1)_{2n-k+1}},
\displaybreak[2]\\
F(n,k)
&=(20n^2-12nk+8n-2k+1)
\frac{(-1)^{n+k}(\frac12)_n^3(\frac12)_{n+k}(\frac12)_{n-k}}
{2^{2n}(1)_n^3(1)_{n-k}^2(\frac12)_k^2},
\displaybreak[2]\\
F(n,k)
&=(3280n^4-4592n^3k+2160n^2k^2-336nk^3+4000n^3-3816n^2k
\\ &\qquad
+1008nk^2-40k^3+1592n^2-884nk+92k^2+232n-62k+13)
\\ &\quad\times 
\frac{(-1)^{n+k}(\frac12)_n^3(\frac12)_{n+k}(\frac12)_{n-k}^3}
{2^{6n}(1)_n^3(1)_{2n-k+1}^2},
\end{align*}
respectively. Nevertheless, following the lines of the above
proofs of Theorems~\ref{t1} and~\ref{t2}, we can show the
validity of~\eqref{A02b} modulo~$p^2$,
of~\eqref{A03b} and~\eqref{A04b} modulo~$p$,
of~\eqref{B01b} modulo~$p^4$, and of~\eqref{B02b} modulo~$p^2$.

We wonder whether there exists a deep general theory
behind all these Ramanujan-type supercongruences,
like the theory of overconvergent $p$-adic modular forms,
which could replace the theory of modular forms used
in the proofs of Ramanujan's formulas for $1/\pi$.

We also wonder whether the methods of the classical theory
of hypergeometric transformations can be applied
to prove the supercongruences mentioned above.
We have at least two successful examples of their application:
D.~McCarthy and R.~Osburn \cite{MO} use them to prove
a different supercongruence conjectured in~\cite{VH1},
and E.~Mortenson~\cite{Mo} gives a similar proof
of the supercongruence~\eqref{A01b}.

\medskip
\textbf{Acknowledgements.}
This work was written during the author's visits in
the Fakult\"at f\"ur Mathematik, Universit\"at Wien (March 2008),
and the Department of Mathematics, National University
of Singapore (May 2008). The author would like to thank Heng Huat Chan
and Christian Krattenthaler for useful conversations on the subject.
The author is also thankful to Frits Beukers, Jonathan Sondow,
Jan Stienstra and Don Zagier for helpful comments,
and to Eric Mortenson for pointing out the references \cite{Mo} and~\cite{VH1}
which became one of the starting points of this project.

\end{document}